\documentclass{article}

\newcommand{\spell}[2]{#1} 

\def\thetitle{Nash's theorem via Günther's trick}

\usepackage{lectures1}
\usepackage[colorlinks=true,
citecolor=black,
linkcolor=black,
anchorcolor=black,
filecolor=black,
menucolor=black,
urlcolor=black,
pdftitle={\thetitle},
pdfsubject={Geometry},
pdfauthor={Anton Petrunin}
]{hyperref}
\makeindex

\spell{}{\makeatletter\usepackage{environ}\RenewEnviron{figure}[1][]{ [PICTURE] }\RenewEnviron{figure*}[1][]{ [PICTURE] }\RenewEnviron{wrapfigure}[1][]{ [PICTURE] }\RenewEnviron{Figure}[1][]{ [PICTURE] }\usepackage{emptypage}\pagestyle{empty}\makeatother\usepackage[none]{hyphenat}\def\emph{\textit}\renewcommand{\footnote}[1]{\textup{\ (#1)}}\renewcommand{\frac}[2]{(#1)/(#2)}\renewcommand{\tfrac}[2]{(#1/#2)}\everydisplay{\textstyle}}

\begin{document}

\title{\thetitle}
\author{Anton Petrunin}
\date{}
\maketitle

\thispagestyle{empty}

\begin{abstract}
We present a proof of Nash's smooth embedding theorem with an emphasis on accessibility and rigor.
\end{abstract}

Riemannian geometry has two fundamental theorems.
One is about the Levi-Civita connection, and the other is a result of John Nash~\cite{nash-1956}, which states that \textit{any Riemannian manifold is isometric to a smooth submanifold of a Euclidean space of large dimension}.
In other words, Riemannian manifolds can be defined as submanifolds of Euclidean spaces with induced metrics.

Nash used a version of Newton's method.
A straightforward application of this method runs into the so-called loss-of-derivatives problem: the correction terms require control of more derivatives than the previous approximation provides.
Nash overcame this by introducing delicate smoothing procedures, unprecedented at the time.
This construction led to the Nash--Moser inverse function theorem, which has found many applications.

Later, Matthias Günther found a magic way that made the loss-of-derivatives problem disappear~\cite{gunter}.
His proof uses a contraction argument, as in the inverse function theorem in Banach spaces.

The proof presented below was tried out on several occasions in lectures and gradually polished along the way.
Except for the inverse function theorem, I mostly follow Nash's proof;
I also use several cosmetic tweaks by
Deane Yang~\cite{yang},
Terence Tao~\cite{tao}, and
Ralph Howard~\cite{howard}.
Günther's approach is also covered in two textbooks: by
Michael Taylor~\cite{taylor} and by
Qing Han and Jia-Xing Hong~\cite{han-hong}.

The reader is assumed to be familiar with smooth manifolds (charts, atlases, coverings, partitions of unity, and tensor fields).
Other results, including the Schauder estimates and the Whitney embedding theorem, can be treated as black boxes.

We formulate the embedding problem in Section~\ref{sec1}.
Section~\ref{sec1+} introduces Nash's twist and solves an approximate version of the problem.
Section~\ref{sec:Reductions} consists of several steps that reduce the embedding problem to the realization of small perturbations of a metric on the torus.
In Section~\ref{sec:free}, we introduce free maps and prove Nash's lemma.
In the following section, we formulate Günther's lemma and use it to realize small perturbations of a metric, thereby proving the embedding theorem.
The final section proves Günther's lemma.
This is the main part of the paper;
the rest is needed to show the power of this lemma.

Here is a list of globally fixed notation:
\spell{\begin{multicols}{3}}{}
\footnotesize
\begin{list}{}{\leftmargin=.5em}
\item $v\oplus w$, \pageref{oplus}
\item $dw$, $w^*$, $\vec x w$, \pageref{dw}
\item $(d\psi)^2$, \pageref{d2}
\item $\Delta u$, \pageref{Du}
\item $f$, \pageref{f}
\item $L$, $L_{ij}$, $M$, \twopageref{LM}{Lij}
\item $N$, \pageref{Nn}
\item $Q$, $Q_{ij}$, $\tilde Q$, \twopageref{Q}{tildeQ}
\item $\T$, $\T^\perp$, $\T^2$, $\N$, \pageref{TTTN}
\end{list}
\spell{\end{multicols}}{}

\section{Formulation}\label{sec1}

\paragraph{Induced metric.}
Given a smooth $n$-dimensional manifold $\Omega$, we denote by $\T\Omega$ and $\T_p\Omega$
its tangent bundle and the tangent space at a point $p \in \Omega$, respectively. \index{4t@$\T\Omega$, $\T_p\Omega$}\label{T}
A \index{metric tensor}\emph{metric tensor}, say $g$, on $\Omega$ is a smooth field of symmetric bilinear forms on $\T\Omega$.
Here and further, \emph{smooth} means~$C^\infty$.

Let $\bar\Omega$ be another smooth manifold, and let $w\colon\Omega \to \bar\Omega$ be a smooth map.
The derivative of $w$ in the direction of a tangent vector $\vec x \in \T_p\Omega$ is denoted by \label{dw}\index{4d@$w^*$, $dw$, $\vec x w$} $\vec x w = dw(\vec x)$;
the right-hand side is the \index{differential}\emph{differential} $dw\colon\T\Omega \to \T\bar\Omega$ of $w$, evaluated at $\vec x$.
Note that $\vec x w \in \T_{w(p)}\bar\Omega$.

Suppose that $\bar\Omega$ is equipped with a metric tensor $\bar g$.
Then the metric tensor $g$ on $\Omega$ defined by
$g(\vec x, \vec y) = \bar g(\vec x w, \vec y w)$
is called the \index{induced metric tensor}\emph{induced} metric tensor, or the \index{pullback}\emph{pullback} of $\bar g$ to $\Omega$.
This relation can also be written as $g = w^*\bar g$.
Equivalently, we may say that $w\colon(\Omega, g)\to(\bar\Omega, \bar g)$ is \index{isometric map}\emph{isometric}.

A metric tensor $g$ is called \index{Riemannian metric tensor}\emph{Riemannian} if it is \textit{positive-definite};
that is, $g(\vec x, \vec x) > 0$ for any nonzero tangent vector $\vec x$.
A smooth manifold equipped with a Riemannian metric is called a \index{Riemannian manifold}\emph{Riemannian manifold}.

Notice that if $g = w^* \bar g$ is a Riemannian metric, then $w$ must be regular; that is, the differential $d_p w$ has rank $n$ at every point $p \in \Omega$.
In particular, $w$ is an immersion.
Conversely, if $w$ is regular and $\bar g$ is Riemannian, then so is $g = w^* \bar g$.

We consider the space $\RR^d$ equipped with the metric tensor defined by the standard scalar product:
\[\bar g(\vec x, \vec y) \df \langle \vec x, \vec y \rangle = x_1 \cdot y_1 + \dots + x_d \cdot y_d,\]
where $x_i$ and $y_i$ denote the coordinates of vectors $\vec x, \vec y \in \RR^d$.
Note that $\bar g$ is Riemannian.

\begin{thm}{Nash's embedding theorem}\label{thm:main}
Any $n$-dimensional Riemannian manifold $(\Omega,g)$ admits a smooth isometric embedding into $\RR^d$ for some $d$.
Moreover, $d$ can be bounded in terms of~$n$.
\end{thm}

We give a complete proof only for compact manifolds and indicate the reduction of the theorem to the compact case;
see Step 4 in the next section.

\paragraph{Remark on \textit{d}.}
I do not make any effort to optimize the bound on $d$ in terms of~$n$.
Tracing the estimates in the argument yields that $d = 10 \cdot n^4$ suffices.

For compact surfaces one can take $d=5$;
this seems to be the only interesting case when the optimal dimension is known.
The existence of isometric embeddings into $\RR^5$ is proved in \cite[3.2.4(D)]{gromov-1986}, while there is no isometric embedding into $\RR^4$ of $\RP^2$ with the standard metric \cite[Appendix 4]{gromov-rokhlin}.

The bound $d\ge \tbinom{n+1}{2}$ is proved in \cite[Appendix 1]{gromov-rokhlin}.
It means that the system
\[
\langle \partial_i w, \partial_j w \rangle = g_{ij},
\eqlbl{eq:<ww>=g}
\]
cannot be overdetermined.
This might be the optimal bound for large $n$ and for the local version of the theorem; that is, \textit{any point in a smooth $n$-dimensional Riemannian manifold has a neighborhood $U$ that admits an isometric embedding into $\RR^d$}.

The system \ref{eq:<ww>=g} makes perfect sense for $C^1$ maps $w$, but the proof of the bound $d\ge \tbinom{n+1}{2}$ uses that $w$ is $C^k$-smooth for sufficiently large $k$.
Without this assumption, the statement does not hold.
Indeed, by the Nash--Kuiper theorem \cite{nash-1954,kuiper-1955} \ref{eq:<ww>=g} has a local $C^1$ solution for $d = n+1$.

\section{Approximate theorem}\label{sec1+}

\paragraph{Sum and direct sum.}
Consider two smooth maps $v_1\colon\Omega\z\to\RR^{d_1}$ and $v_2\colon\Omega\z\to\RR^{d_2}$.
The map \label{oplus}\index{$v_1\oplus v_2$}$v_1\oplus v_2\colon p\mapsto (v_1(p),v_2(p))$ will be called the \index{direct sum}\emph{direct sum} of $v_1$ and $v_2$;
it maps $\Omega$ to $\RR^{d_1+d_2}\z=\RR^{d_1}\oplus\RR^{d_2}$.

Suppose that metric tensors $g$, $g_1$, and $g_2$ are induced by $v_1\oplus v_2$, $v_1$, and $v_2$, respectively.
Then
\[g=g_1+g_2.\]

Furthermore, if $d_1=d_2$, then we can consider the usual sum $v_1+v_2$.
Recall that the \index{support}\emph{support} of a map is the closure of the set where it takes nonzero values.
Now suppose that $g$, $g_1$, and $g_2$ are induced by $v_1+ v_2$, $v_1$, and $v_2$, and the maps $v_1$ and $v_2$ have disjoint supports.
Then again we have
\[g=g_1+g_2.\]

These two observations will be used repeatedly without notice.

\paragraph{Nash's twist.}
Let $\phi$ and $\psi$ be smooth functions on a smooth manifold $\Omega$.
Given $r>0$, denote by $\SSS^1_r$ the circle of radius $r$ in $\RR^2$.

{

\begin{wrapfigure}[9]{r}{30mm}
\vskip-0mm
\centering
\includegraphics{mppics/pic-15}
\end{wrapfigure}

\mbox{\index{Nash's twist}\emph{Nash's twist}} $\Theta\colon\Omega \to \RR^2$ of the triple $(r,\phi,\psi)$ is defined as the composition
\[\Omega\xrightarrow{\psi}\RR\xrightarrow{\ell_r}\SSS^1_r\xrightarrow{\times \phi}\RR^2,\]
where $\ell_r\colon\RR\to\SSS^1_r$ is a length-preserving covering map, say $\ell_r(x)=(r\cdot\cos\tfrac xr,r\cdot\sin\tfrac xr)$,
and $\times \phi$ denotes multiplication by $\phi$; so
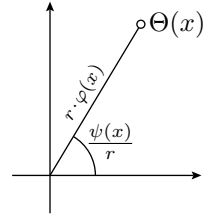
\[\Theta(x)
\df
\phi(x)\cdot (\ell_r\circ \psi(x)).\]

}

\begin{thm}{Claim}\label{clm:twist}
Nash's twist $\Theta$ for the triple $(r,\phi,\psi)$ induces the following metric
\[g=\phi^2\cdot (d\psi)^2+r^2\cdot(d\phi)^2.\]
Here \label{d2}\index{4d@$(d\psi)^2$}$(d\psi)^2$ denotes the metric tensor induced by the function
$\psi\colon\Omega \z\to \RR$; that is,
\[(d\psi)^2(\vec x,\vec y)\df d\psi(\vec x)\cdot d \psi(\vec y)=(\vec x\psi)\cdot (\vec y\psi).\]
\end{thm}

\parit{Computations.}
Suppose that $\vec x$ is a tangent vector on $\Omega$.
Then
\begin{align*}
\vec x\Theta&=\vec x(\phi\cdot \ell_r\circ \psi)=
\\
&=(\vec x\phi)\cdot(\ell_r\circ \psi) + \phi\cdot (\ell_r'\circ \psi)\cdot(\vec x\psi).
\end{align*}

Observe that $|\ell_r|=r$, $|\ell'_r|=1$, and $\ell_r\perp \ell'_r$.
Therefore,
\[g(\vec x,\vec y)=\langle \vec x\Theta,\vec y\Theta\rangle =r^2\cdot (\vec x\phi)\cdot (\vec y\phi)+\phi^2\cdot (\vec x\psi)\cdot (\vec y\psi).
\eqno\qed\]

\paragraph{Approximate version.}
The following statement can be regarded as an approximate version of the main theorem (\ref{thm:main}).
It implies in particular that the set of metrics induced by smooth maps into Euclidean space is dense in the space of all Riemannian metrics.

\begin{thm}{Proposition}\label{prop:approx-nash}
Let $g$ be a Riemannian metric on a smooth $n$-manifold $\Omega$.
Then there exists a metric tensor $h$ on $\Omega$ and a one-parameter family of smooth maps $w_t\colon\Omega \to \RR^d$ for $t>0$ with the induced metrics $g_t=g+t\cdot h$.

Moreover, $d$ depends only on $n$.
\end{thm}

The proof will produce a nonnegative metric tensor $h$, but this is not going to be used.

Observe that $g_t\to g$ as $t \to 0$.
Therefore, one might be tempted to take the limit of $w_t$ as $t \to 0$.
However, as we will see, the maps $w_t$ constructed in the proof converge to a constant;
thus, the limit does not solve the problem.

Let $w\colon\Omega \to \RR^d$ be a smooth map.
Suppose that $w_1,\dots,w_d$ are the coordinate functions of $w$.
Then the induced metric tensor $g$ can be written as
\[g=(dw_1)^2+\dots+(dw_d)^2.\eqlbl{eq:nash=dw^2}\]
Solving \ref{eq:nash=dw^2} is equivalent to the embedding problem.
The following weaker form of \ref{eq:nash=dw^2} will play a key role in the proof of \ref{prop:approx-nash}.

\begin{thm}{Proposition}\label{prop:phi-f}
Let $(\Omega,g)$ be an $n$-dimensional Riemannian manifold.
Then there are smooth functions
$\phi_1,\z\dots,\phi_d,\psi_1,\z\dots,\psi_d$ on $\Omega$
such that
\[g=\phi_1^2 \cdot (d\psi_1)^2+\dots+\phi_d^2\cdot (d\psi_d)^2.\eqlbl{eq:g=phidf}\]

Moreover, $d$ depends only on $n$.
\end{thm}

\paragraph{Proof.}
The metric tensor $g$ can be written in local coordinates $(x_1,\z\dots,x_n)$ as
\[
g = \sum_{i, j} g_{ij} \cdot dx_i \odot dx_j,
\eqlbl{g=sum}
\]
where $g_{ij} = g_{ji}$ are smooth functions of $(x_1,\dots,x_n)$ and $\odot$\label{odot} denotes the symmetric tensor product.

Since $g$ is Riemannian, at any point $p \in \Omega$, we can choose a chart so that the vectors $\partial_i$ are orthonormal at $p$;
that is, $g_{ii} = 1$ and $g_{ij} = 0$ at $p$ for all $i \ne j$.
Since $g_{ij}$ are smooth functions, for any $\eps > 0$, we can find a neighborhood $U \ni p$ such that
\[
g_{ii} \lessgtr 1 \pm \eps
\quad\text{and}\quad
g_{ij}\lessgtr \pm \eps
\eqlbl{alpha=1}\]
for all $i\ne j$ and any point in $U$.

Observe that
\[\pm dx_i\odot dx_j=\tfrac12\cdot (dx_i\pm dx_j)^2-\tfrac12\cdot(dx_i)^2-\tfrac12\cdot(dx_j)^2.\]
Let us plug this formula into \ref{g=sum} and take \ref{alpha=1} into account.
Assuming that $\eps$ is small, we will get that $g$ is a linear combination with positive coefficients of the metric tensors
$(dx_i)^2$ and $(dx_i\pm dx_j)^2$.
In other words, we can take the functions $x_i$ and $x_i\pm x_j$ for all $i< j$ as our functions $\psi_1,\dots,\psi_d$ (each should be extendable to a smooth function on the whole manifold)
and find $\phi_1,\dots,\phi_d$ such that \ref{eq:g=phidf} holds in a neighborhood of~$p$.
Applying a partition of unity, we get the global statement.

It suffices to take $n^2$ functions in each chart.
So, if $\Omega$ is covered by $m$ charts, then we may take $d = n^2 \cdot m$.
Moreover, the same bound applies if the charts can be colored using $m$ colors such that distinct charts of the same color do not overlap.
By the colored version of Lebesgue covering dimension \cite[Theorem 2]{ostrand},
$m=n+1$ colors are sufficient;
therefore, we can take $d \z= n^2 \cdot (n+1)$.
\qeds

\paragraph{Proof of \ref{prop:approx-nash}.}
Suppose that $\phi_1,\z\dots,\phi_d,\psi_1,\z\dots,\psi_d$ are provided by \ref{prop:phi-f}.
Denote by $\Theta_i$ Nash's twist for the triple $(\sqrt{t},\phi_i,\psi_i)$.
By \ref{clm:twist},
\[w_t=\Theta_1\oplus\dots\oplus \Theta_d\]
is the needed family of maps with
$h=(d\phi_1)^2+\dots+(d\phi_d)^2$.
\qedsf

\paragraph{A pseudoeuclidean digression.}
Let us denote by $\RR^r\ominus\RR^s$\index{3c@$\RR^r\ominus\RR^s$, $v\ominus w$} the pseudoeuclidean space with signature $(r,s)$;
that is, the space $\RR^{r+s}$ with the metric tensor
\[\bar g(\vec x,\vec y)=\langle \vec x,\vec y\rangle=x_1\cdot y_1+\dots+x_r\cdot y_r-x_{r+1}\cdot y_{r+1}-\dots-x_{r+s}\cdot y_{r+s},\]
where $x_i$ and $y_i$ denote the coordinates of vectors $\vec x$ and $\vec y$ in $\RR^{r+s}$.

\begin{thm}{Problem}
Show that any metric tensor $g$ on a smooth manifold $\Omega$ is induced by a smooth embedding $\Omega\hookrightarrow \RR^r\ominus\RR^s$ for large $r$ and $s$.
\end{thm}

\parit{Solution based on Nash's theorem.}
The metric tensor can be written as the difference $g = g_1 - g_2$ of two Riemannian metrics $g_1$ and $g_2$ on $\Omega$.
By Nash's theorem, we can find two smooth embeddings $w_1, w_2 \colon \Omega \z\hookrightarrow \RR^d$ whose induced metrics are $g_1$ and $g_2$, respectively.
It remains to observe that $g$ is induced by \label{ominus}
\[w_1\ominus w_2\colon p\mapsto(w_1(p),w_2(p))\in \RR^d\ominus\RR^d.\eqno{\qed}\]

The above solution uses the main theorem, which we have not yet proved.
The following solution is more direct and relies only on Nash's twist.

\parit{Direct solution.}
Assume $g$ is Riemannian, and let $\phi_1,\z\dots,\phi_d,\psi_1,\z\dots,\psi_d$ be provided by \ref{prop:phi-f}.
Let $\Theta_i$ be Nash's twist for the triple $(1,\phi_i,\psi_i)$, and
\begin{align*}
h&=(d\phi_1)^2+\dots+(d\phi_d)^2,
\\
\Theta&=\Theta_1\oplus\dots\oplus \Theta_d\colon\Omega\to \RR^{2\cdot d},
\\
\phi&=\phi_1\oplus \dots\oplus \phi_d\colon\Omega\to \RR^{d}.
\end{align*}
According to \ref{clm:twist}, the map $\Theta$ induces $g+h$.
Note that $\phi$ induces $h$.
Therefore $w=\Theta\ominus \phi\colon\Omega\to \RR^{2\cdot d}\ominus\RR^d$ induces $g$.

To handle the general case, express the metric as a difference of two Riemannian metrics: $g = g_1 - g_2$.
As shown above, there exist smooth maps $w_1$ and $w_2$ into pseudoeuclidean spaces inducing $g_1$ and $g_2$, respectively.
It follows that $w_1 \ominus w_2$ induces $g$.

Finally, choose $w_3\colon\Omega\hookrightarrow \RR^{2\cdot n+1}$ provided by the Whitney embedding theorem, and observe that $(w_1 \oplus w_3)\ominus (w_2\oplus w_3)$ is the required smooth isometric embedding.
\qeds

\section{Reductions}\label{sec:Reductions}

\paragraph{Realization of small perturbations of a metric.}
Let us denote by $\TT^n=\RR^n/\ZZ^n$ the $n$-dimensional torus.

\begin{thm}{Theorem}\label{thm:restricted}
There exists a Riemannian metric $g_0$ on $\TT^n$ such that
any metric $g$ that is sufficiently close to $g_0$ in the $C^\infty$ topology is induced by a smooth map $\TT^n\to \RR^d$ for some $d$ that depends only on~$n$.
\end{thm}

This theorem will be proved in the following sections;
note that this is a very special case of the main theorem (\ref{thm:main}).
Now we will reduce the main theorem to this restricted version.

\paragraph{Step 1.}
Let us use \ref{prop:approx-nash} to construct an isometric immersion for any metric on the torus; namely,
\begin{clm}{}
given a Riemannian metric $g$ on $\TT^n$, there exists an isometric immersion $(\TT^n, g) \looparrowright \RR^d$ for some $d$ depending only on~$n$.
\end{clm}

Indeed, let $g_0$ be the metric from \ref{thm:restricted}.
After rescaling if necessary, we may assume that $g_0<g$; that is, $g_0(\vec x,\vec x)<g(\vec x,\vec x)$ for any nonzero tangent vector $\vec x$.
In particular, the difference $g_1=g-g_0$ is a Riemannian metric.

By \ref{prop:approx-nash}, there exists a metric tensor $h$ such that, given $t>0$, there is a smooth map $w_1\colon\TT^n\to \RR^{d_1}$ that induces $g_1+t\cdot h$.

By \ref{thm:restricted}, for any small $t>0$, there is a smooth map $w_0\colon\TT^n\to \RR^{d_0}$ that induces the metric $g_0-t \cdot h$.
Therefore, $w_0\oplus w_1$ induces
\[g=g_0-t\cdot h+g_1+t\cdot h.\]

\paragraph{Step 2.}
Now let us improve the statement to an embedding; namely,
\begin{clm}{}
given a Riemannian metric $g$ on $\TT^n$, there is an isometric embedding $(\TT^n,g)\hookrightarrow\RR^d$ for some $d$ depending only on~$n$.
\end{clm}

Choose a smooth regular embedding $w_1\colon\TT^n\hookrightarrow \RR^{n+1}$.
Let $g_1$ be the metric induced by $w_1$.
Since $w_1$ is regular, $g_1$ must be Riemannian.

We can assume that $g>g_1$; in other words, the metric $g_2=g-g_1$ is Riemannian.
If this is not the case, replace $w_1$ with its rescaling $\varepsilon \cdot w_1$ for small $\varepsilon > 0$.

Step 1 provides an isometric immersion $w_2\colon(\TT^n,g_2)\z\looparrowright \RR^d$.
Since $g=g_1+g_2$,
\[w=w_1\oplus w_2\colon (\TT^n,g)\to \RR^{d+ n+1}\]
is an isometric embedding; $w$ is an embedding since $w_1$ is.

\paragraph{Step 3.} Now we extend the statement to all compact manifolds; namely,
\begin{clm}{}
any compact $n$-dimensional Riemannian manifold $(\Omega,g)$ admits a smooth isometric embedding into $\RR^d$ for some $d$ that depends only on~$n$.
\end{clm}

The Whitney theorem gives a smooth embedding $\iota\colon\Omega\hookrightarrow\TT^{2\cdot n+1}$.

The Riemannian metric on $\Omega$ can be extended to a Riemannian metric on the torus.
This follows, for example, by choosing a metric on a tubular neighborhood that restricts to the given one on \(\Omega\), and then patching it with the standard metric on the torus via a partition of unity.
It remains to take a smooth isometric embedding $w\colon (\TT^{2\cdot n+1}, \bar g) \z\hookrightarrow \RR^d$ provided by Step 2 and observe that the composition $w \circ \iota\colon(\Omega, g) \z\hookrightarrow \RR^d$ is an isometric embedding.

\paragraph{Step 4.}
Now we sketch the reduction of the main theorem to the compact case.
A complete proof is given in \cite[Part~D]{nash-1956}.

Applying an argument similar to Step 2, we can reduce the question to finding an isometric immersion of $(\Omega,g)$.

Further, we need to find a sequence of maps $r_{\mathfrak{i}}\colon\Omega \to \SSS^n$ and a sequence of Riemannian metrics $h_{\mathfrak{i}}$ on $\SSS^n$ such that
\[g=\sum_{\mathfrak{i}} r_{\mathfrak{i}}^*h_{\mathfrak{i}}\]
and $r_{\mathfrak{i}}$ is constant, say $c_{\mathfrak{i}}$, outside of a compact subset $K_{\mathfrak{i}}\subset \Omega$.
This part is done by applying a partition of unity.

In addition, we have to arrange this construction so that the index set $\mathfrak{I}$ can be colored with $n+1$ colors such that $K_{\mathfrak{i}}\cap K_{\mathfrak{j}}=\emptyset$ for distinct indices $\mathfrak{i}$ and $\mathfrak{j}$ of the same color.
More precisely, $\mathfrak{I}$ can be partitioned into subsets $\mathfrak{I}_0, \dots, \mathfrak{I}_n$ such that for any $k$ and any distinct $\mathfrak{i}, \mathfrak{j} \z\in \mathfrak{I}_k$, we have
$K_{\mathfrak{i}} \cap K_{\mathfrak{j}} = \emptyset$.
This is very similar to the colored version of Lebesgue covering dimension.

Applying the compact case of the theorem to each $(\SSS^n,h_{\mathfrak{i}})$, we get a sequence of isometric embeddings $w_{\mathfrak{i}}\colon(\SSS^n,h_{\mathfrak{i}})\hookrightarrow \RR^d$.
By shifting $w_{\mathfrak{i}}$, we can assume that $w_{\mathfrak{i}}(c_{\mathfrak{i}})=0$, and therefore $w_{\mathfrak{i}}\circ r_{\mathfrak{i}}$ has support in $K_{\mathfrak{i}}$.
Now note that the map
\[\biggl(\,\sum_{\mathfrak{i}\in \mathfrak{I}_0}w_{\mathfrak{i}}\circ r_{\mathfrak{i}}\biggr)\oplus\dots\oplus
\biggl(\,\sum_{\mathfrak{i}\in \mathfrak{I}_n}w_{\mathfrak{i}}\circ r_{\mathfrak{i}}\biggr)\]
gives the required isometric immersion $(\Omega,g)\looparrowright \RR^{(n+1)\cdot d}$.

\section{Free maps}\label{sec:free}

From now on, we look only at the torus $\TT^n\z=\RR^n/\ZZ^n$, which is equipped with cyclic coordinates and globally defined partial derivatives $\partial_1,\dots,\partial_n$.
Furthermore, a metric tensor on $\TT^n$ can be (and will) described by a smooth map $g\colon\TT^n\to\RR^N$ with components $g_{ij}$ for $1\le i\le j\le n$; here and further \label{Nn}
\[N\df\tfrac{n\cdot(n+1)}2.\]

\paragraph{Free maps.}
Recall that a smooth map $f\colon\TT^n \to \RR^d$ is \index{regular map}\emph{regular} if its differential $df$ has rank $n$ at every point.
In other words, for any point $p\in \TT^n$, the first-order partial derivatives $\partial_1 f,\dots,\partial_n f$ are linearly independent.

A map $f\colon\TT^n \to \RR^d$ is called \index{free map}\emph{free} if an analogous property holds for both first- and second-order partial derivatives:
that is, at any point $p\in \TT^n$, all $n + N$ partial derivatives
$\partial_i f$,
$\partial_{ij}f$ for $i\le j$
are linearly independent.

Let $\psi_\ell\colon\TT^n\to\SSS^1=\RR/\ZZ$ be a map induced by a linear function $\ell\colon\RR^n\to\RR$.
Let us regard $\SSS^1$ as a subset of $\RR^2$ and consider the map $f\colon\TT^n\z\to \RR^{n\cdot(n+1)}$
defined as the direct sum of the maps
$\psi_{x_i}$ and $\psi_{x_i+x_j}$, where $x_i$ are the coordinate functions on $\RR^n$ and $i<j$.
It is straightforward to check that $f$ is free, so we have proved the following claim.

\begin{thm}{Claim}\label{ex:free-composition}
There is a smooth free embedding $f\colon\TT^n\hookrightarrow \RR^{n\cdot(n+1)}$.
\end{thm}

From now on a free map $f\colon\TT^n\hookrightarrow \RR^d$\label{f} is fixed and we may treat $\TT^n$ as a submanifold of $\RR^d$.
The tangent space $\T_p=\T_p\TT^n$ is the $n$-dimensional space spanned by the partial derivatives $\partial_1f,\z\dots, \partial_nf$ at~$p$.
Let $\T_p^\perp$ denote the orthogonal complement of $\T_p$ in~$\RR^d$.
Further, define the osculating space $\T^2_p$ as the span of all first- and second-order partial derivatives $\partial_i f$ and $\partial_{ij} f$ at $p$, and let us set $\N_p \z\df \T^2_p \cap \T_p^\perp$.

Since $f$ is free, the dimensions of the spaces $\T_p$, $\T_p^\perp$, $\T^2_p$, and $\N_p$ are $n$, $d - n$, $n + N$, and $N$, respectively.
We will write $\T$, $\T^\perp$, $\T^2$, and $\N$ for the corresponding vector bundles over $\TT^n$.
\index{4t@$\T$, $\T^\perp$, $\T^2$, $\N$}\label{TTTN}

\paragraph{Hölder norms.}
Fix some $\alpha\in(0,1)$ once and forever, say $\alpha=\tfrac12$.
Denote by $|x-y|_{\TT^n}$ the standard distance between points $x,y\z\in\TT^n$.
Let us recall the definition of the Hölder norms on an open subset $W\subset \TT^n$: \label{C(T)}
\begin{align*}
\|u\|_{C^{0,\alpha}(W)}&\df\sup_{x}\{\,|u(x)|\,\}+\sup_{x\ne y}\left\{\,\frac{|u(x)-u(y)|}{|x-y|_{\TT^n}^\alpha}\,\right\},
\\
\|u\|_{C^{k,\alpha}(W)}&\df\max_{|\mathcal{I}|\le k}\left\{\,\|\partial^{\mathcal{I}}u\|_{C^{0,\alpha}(W)}\,\right\},
\end{align*}
where $x,y\in W$, and $\mathcal{I} = (i_1,\dots,i_n)$ denotes a multi-index; so, $\partial^{\mathcal{I}}$ is the partial derivative $\partial_1^{i_1}\dots\partial_n^{i_n}$ of order $|\mathcal{I}| \z\df i_1 +\dots+i_n$.
The space of functions $W\to \RR$ with finite ${C^{k,\alpha}(W)}$-norm is called the $(k,\alpha)$-Hölder space on $W$ and is denoted by $C^{k,\alpha}(W)$; \index{$C^{k,\alpha}$, $\Arrowvert u\Arrowvert_{C^{k,\alpha}}$} we may suppress $W$ in this notation.

Let us define $C^{k,\alpha}(\TT^n,\RR^d)$ as the space of maps $u\colon\TT^n\z\to \RR^d$ whose coordinate functions $u_1,\dots, u_d$ belong to $C^{k,\alpha}(\TT^n)$,
and set
\[\|u\|_{C^{k,\alpha}}\df\max \{\,\|u_1\|_{C^{k,\alpha}},\dots,\|u_d\|_{C^{k,\alpha}} \,\}.\]

\paragraph{Q-form.}
We will need a somewhat fancy notation for the induced metric.
Let us define a symmetric bilinear form $Q$\index{1Q@$Q(v,w)$, $Q_{ij}(v,w)$}\label{Q} that takes two smooth maps $v, w\colon\TT^n \to \RR^d$ and returns the following metric tensor on~$\TT^n$:
\[
g(\vec x, \vec y)
\df
\tfrac{1}{2} \cdot [\langle \vec x v, \vec y w \rangle + \langle \vec x w, \vec y v \rangle].
\]
Recall that the cyclic coordinates on $\TT^n$ allow us to write a metric tensor as a smooth map $\TT^n\to\RR^N$.
The components of the metric tensor can be written as
\[
g_{ij} = Q_{ij}(v, w) = \tfrac{1}{2} \cdot [\langle \partial_i v, \partial_j w \rangle+\langle \partial_i w, \partial_j v \rangle].
\]

Notice that a metric tensor $g$ is induced by a smooth map $w\colon\TT^n \z\to \RR^d$
if and only if $g = Q(w, w)$.

\paragraph{Linearization.}
Recall that $f\colon\TT^n \hookrightarrow\RR^d$ is a fixed smooth free embedding
(which exists by \ref{ex:free-composition}).
Let us define the linear operator
\[L\colon C^\infty(\TT^n,\RR^d)\to C^\infty(\TT^n,\RR^N)\qquad
\text{by}
\qquad L(w)\df 2\cdot Q(f,w).\]

We will see soon that realization of small perturbations of a metric follows if the quadratic equation
\[L(w)+Q(w,w)=h\]
has a solution for small metric tensors $h$.
The following lemma provides a solution of the corresponding linear equation $L(w)=h$.

\begin{thm}{Nash's lemma}\label{lem:nash-guenter:nash}
The linear operator $L$ admits a right inverse \index{1Q@$L(w)$, $M(h)$}\label{LM}
\[M\colon C^\infty(\TT^n,\RR^N)\to C^\infty(\TT^n,\RR^d)\]
such that
\[\|Mh\|_{C^{2,\alpha}}\le a\cdot \|h\|_{C^{2,\alpha}}\]
for some constant $a$ and any $h\in C^\infty(\TT^n,\RR^N)$.

Moreover, $Mh$ is an $\N$-field;
that is, $Mh(p)\in \N_p=\T_p^\perp\cap \T_p^2$ for all $p\in \TT^n$.
\end{thm}

\paragraph{Proof.}
Choose a metric tensor $h$ on $\TT^n$.
Since all $\partial_i f$,
$\partial_{ij}f$
are linearly independent,
the equations
$\langle\partial_{ij}f,w\rangle=-\tfrac12\cdot h_{ij}$
uniquely define an $\N$-field $w$.
Equivalently, $w$ can be defined as a linear combination of $\partial_i f$,
$\partial_{ij}f$ such that
\[
\begin{aligned}
\langle\partial_if,w\rangle&=0,
&
-2\cdot\langle\partial_{ij}f,w\rangle&=h_{ij}
\end{aligned}
\eqlbl{eq:w(h)}
\]
at each point.

Notice that \label{Lij}
\begin{align*}
L_{ij}(w)&=2\cdot Q_{ij}(f,w)=
\\
&=\langle\partial_if,\partial_jw\rangle+\langle\partial_jf,\partial_iw\rangle=
\\
&=\partial_j\langle\partial_if,w\rangle+\partial_i\langle\partial_jf,w\rangle -2\cdot \langle\partial_{ij}f,w\rangle.
\end{align*}
By \ref{eq:w(h)}, the operator $M\colon h \mapsto w$ serves as the desired right inverse.
\qeds

\section{Realization of small perturbations of a metric}\label{sec:perturbation}

\begin{thm}{Günther's lemma}\label{lem:nash-guenter:guenter}
There exists a symmetric bilinear form \index{1Q@$\tilde Q(v,w)$}\label{tildeQ}
\[\tilde Q\colon C^\infty(\TT^n,\RR^d)\times C^\infty(\TT^n,\RR^d) \to C^\infty(\TT^n,\RR^d)\]
such that $L\tilde Q=Q$ and for a sequence of constants $b_2,b_3,\dots$ we have
\begin{align*}
\|\tilde Q(w,w)\|_{C^{2,\alpha}}&\le b_2\cdot \|w\|_{C^{2,\alpha}}^2,
\\
\|\tilde Q(w,w)\|_{C^{k,\alpha}}
&\le b_2\cdot \|w\|_{C^{2,\alpha}}\cdot \|w\|_{C^{k,\alpha}}
+
b_k\cdot \|w\|_{C^{k-1,\alpha}}^2
\end{align*}
for any $w\in C^\infty(\TT^n,\RR^d)$ and $k\ge 3$.
\end{thm}

The proof is given in the next section.

\paragraph{Realization of small perturbations of a metric.}
We now complete the proof of Nash's theorem, modulo Günther's lemma.
Recall that $f\colon\TT^n \hookrightarrow\RR^d$ is a fixed smooth free embedding and $L(w)\z\df 2\cdot Q(f,w)$.

\begin{thm}{Lemma}\label{lem:perturb}
There exists $r>0$ such that the equation
\[L(w)+Q(w,w)=h\eqlbl{eq:L+Q=h}\]
admits a smooth solution $w\colon\TT^n\to\RR^d$ for any smooth metric tensor~$h$ (not necessarily Riemannian) with $\|h\|_{C^{2,\alpha}}\z<r$.
Moreover, we can assume that
\[\|w\|_{C^{2,\alpha}}\le c\cdot \|h\|_{C^{2,\alpha}}\]
for some constant $c$.
\end{thm}

In Section \ref{sec:Reductions}, we showed that realization of small perturbations of a metric (\ref{thm:restricted}) implies the main theorem (\ref{thm:main}).
Now we will show that the lemma implies realization of small perturbations (\ref{thm:restricted}), and therefore the main theorem (\ref{thm:main}).

\paragraph{Proof of \ref{thm:restricted} modulo \ref{lem:perturb}.}
Recall that $f\colon\TT^n\to\RR^d$ is a fixed free map.
Let $g_0$ be the metric tensor induced by $f$, so $g_0=Q(f,f)$.
Further, let $h=g-g_0$,
so $h$ is a small metric tensor; in particular we can assume that $\|h\|_{C^{2,\alpha}}\z<r$.

It is sufficient to find a smooth map $w\colon\TT^n\to\RR^d$ such that
\[
Q(f+w, f+w) = g_0 + h,
\]
and this equation is equivalent to \ref{eq:L+Q=h}.
Indeed,
\begin{align*}
Q(f+w, f+w)
&= Q(f, f) + 2 \cdot Q(f, w) + Q(w, w)
\\
&= g_0 + L(w) + Q(w, w).
\end{align*}
Hence, \ref{thm:restricted} follows. \qeds

The next proof follows the argument of the inverse function theorem in Banach spaces; see, for example,~\cite{cartan}.

\paragraph{Proof of \ref{lem:perturb} modulo \ref{lem:nash-guenter:guenter}.}
Recall that $M$, $\tilde Q$, $a$, and $b_2$, $b_3,\dots$ are provided by \ref{lem:nash-guenter:nash} and \ref{lem:nash-guenter:guenter}.
Set $R=\tfrac1{10\cdot b_2}$ and $r=\tfrac R{10\cdot a}$.
Assume $\|h\|_{C^{2,\alpha}}\z\le r$.
By \ref{lem:nash-guenter:guenter},
\[\Phi\colon w\mapsto Mh-\tilde Q(w,w)\]
is a contraction on the ball $\cBall[0,R]\subset {C^{2,\alpha}(\TT^n,\RR^d)}$.

Consider the sequence of maps $w_0,w_1,\dots\colon\TT^n\to\RR^d$ defined by $w_0=0$ and $w_{n+1}=\Phi(w_n)$ for all $n$.
Notice that each $w_n$ is smooth.
Since $\Phi$ is a contraction, the sequence $w_n$ converges in $C^{2,\alpha}(\TT^n,\RR^d)$;
denote its limit by $w$.
Then
\[\|w\|_{C^{2,\alpha}}\le R
\qquad\text{and}\qquad
w= Mh-\tilde Q(w,w).
\eqlbl{eq:w=Mh-Q(w,w)}\]

Since $LM=\id$ and $L\tilde Q=Q$, applying $L$ to both sides of the equality \ref{eq:w=Mh-Q(w,w)} yields that $w$ solves \ref{eq:L+Q=h}.

Let us show that $w$ is smooth using the following claim:
\begin{clm}{}\label{clm:wn-bounded}
The sequence $\|w_n\|_{C^{k,\alpha}}$ is bounded for every integer $k\ge 2$.
\end{clm}
Indeed, since the embedding
$C^{k,\alpha}(\TT^n,\RR^d)\hookrightarrow C^k(\TT^n,\RR^d)$ is compact,
a subsequence of $w_n$ converges in $C^k(\TT^n,\RR^d)$;
but its limit must coincide with $w$.
So $w$ is $C^k$-smooth for every $k$ and therefore smooth.

We prove \ref{clm:wn-bounded} by induction on $k$.
The base case $k=2$ has already been proved.
By the induction hypothesis, $\|w_n\|_{C^{k-1,\alpha}}$ is bounded.
In particular, there is a constant $c_k$ such that
\[\|Mh\|_{C^{k,\alpha}}+b_k\cdot \|w_n\|_{C^{k-1,\alpha}}^2\le c_k\]
for all $n$; here $b_k$ is the constant from \ref{lem:nash-guenter:guenter}.

Since $\|w_n\|_{C^{2,\alpha}}\le R=\tfrac1{10\cdot b_2}$, by the second estimate in \ref{lem:nash-guenter:guenter}, we get
\begin{align*}
\|w_{n+1}\|_{C^{k,\alpha}}
&\le
\|Mh\|_{C^{k,\alpha}}
+
b_2\cdot \|w_n\|_{C^{2,\alpha}}\cdot \|w_n\|_{C^{k,\alpha}}
+
b_k\cdot \|w_n\|_{C^{k-1,\alpha}}^2\le
\\
&\le c_k +\tfrac1{10}\cdot\|w_n\|_{C^{k,\alpha}}.
\end{align*}
In particular, $\|w_n\|_{C^{k,\alpha}}\le 2\cdot c_k
\ \Rightarrow\
\|w_{n+1}\|_{C^{k,\alpha}}\le 2\cdot c_k$
for all $n$.
Since $w_0=0$, we get \ref{clm:wn-bounded}.
\qeds

\section{Proof of Günther's lemma}\label{sec:proof(gunther)}

Recall that $L$ and $M$ are defined on page \pageref{LM}.
The form $\tilde Q$ in Günther's lemma can be defined by the formula
\[\tilde Q(w,w)=M(Q(w,w)+L\vec x)-\vec x,\]
where $\vec x$ is a tangent vector field on $\TT^n$ that depends quadratically on $w$;
see the remark right after the proof below.
Since $M$ is a right inverse of $L$, we have $L\tilde Q=Q$.
It remains to choose $\vec x$ so that the inequalities in the lemma hold; see \ref{eq:Delta}, \ref{eq:AB=<C}, and \ref{eq:tildeQ} below.

The correction term $w$ in the proof of \ref{lem:perturb} satisfies the equation
\[w= Mh-\tilde Q(w,w)=M(h-Q(w,w)-L\vec x)+\vec x.\]
By Nash's lemma, $\vec x$ is the tangential part of $w$.
This tangential term is the main departure from Nash's original construction.
Notice that perturbing an embedding in a tangential direction almost does not change the Riemannian metric: to first order, it only changes the parametrization.
For this reason, introducing a tangential term might look strange, but it works.
The same idea works in the DeTurck trick for the Ricci flow \cite[2.2]{brendle} and in Kuiper's improvement of Nash's $C^1$-embedding theorem~\cite{kuiper-1955}.

\paragraph{Smoothing operator.}
Recall that $\alpha\in(0,1)$ is fixed, say $\alpha=\tfrac12$.
Let $u\colon \TT^n\to\RR$ be a smooth function;
its Laplacian is defined by \index{4d@$\Delta u$}\label{Du}
\[\Delta u=\sum_s\partial_{ss}u,\]
where $\partial_1,\dots,\partial_n$ are the standard partial derivatives on $\TT^n$.

\begin{thm}{Proposition}\label{prop:Delta-1}
The operator $(\Delta-1)\colon u\mapsto \Delta u-u$ on $C^\infty(\TT^n)$ has an inverse $(\Delta-1)^{-1}$.
Moreover, there is a constant $c>0$ such that
\[\|(\Delta-1)^{-1}v\|_{C^{k,\alpha}(\TT^n)}\le c\cdot \|v\|_{C^{k-2,\alpha}(\TT^n)}\]
for every integer $k\ge2$ and every function $v\in C^\infty(\TT^n)$.
\end{thm}

\parit{Sketch.}
The existence of $(\Delta-1)^{-1}$ follows from Fourier analysis.
In other words, the equation $\Delta u - u = v$ has a unique smooth solution $u$ for a given $v\in C^\infty(\TT^n)$.
It remains to prove the following inequality with some constant $c$:
\[\|u\|_{C^{k,\alpha}(\TT^n)}
\le
c\cdot\|v\|_{C^{k-2,\alpha}(\TT^n)}.
\eqlbl{eq:biLip}\]

By the maximum principle, we have
\[\|u\|_{L^\infty(\TT^n)}\le\|v\|_{L^\infty(\TT^n)},\eqlbl{eq:maximum-principle}\]
where $\|\ \|_{L^\infty(W)}$ denotes the sup-norm on $W\subset \TT^n$.

Choose a $\tfrac1{10}$-ball $B$ in $\TT^n$ (with the standard metric);
denote by $2\cdot B$ the concentric ball with twice the radius.
Since $2\cdot B$ is isometric to a Euclidean ball, Schauder's estimates can be applied to the restriction $u|_{2\cdot B}$ (see, for example \cite[Lemma 6.1(a)]{gilbarg-trudinger}).
By \ref{eq:maximum-principle}, this implies
\begin{align*}
\|u\|_{C^{2,\alpha}(B)}
&\le c\cdot(\|v\|_{C^{0,\alpha}(2\cdot B)} + \|u\|_{L^\infty(2\cdot B)})\le
 2\cdot c\cdot\|v\|_{C^{0,\alpha}(\TT^n)}.
\end{align*}
Averaging this inequality for $\tfrac1{10}$-balls in $\TT^n$ leads to \ref{eq:biLip} for $k=2$.

Since $\partial_i\Delta=\Delta\partial_i$ for each $i$,
applying the already proved case to the partial derivatives $\partial^{\mathcal{I}}u$ and $\partial^{\mathcal{I}}v$, we obtain \ref{eq:biLip} for all $k$.
\qeds

Notice that $(\Delta-1)^{-1}$ is a smoothing operator.
Moreover, $(\Delta-1)^{-1}$ commutes with partial derivatives; that is,
\[(\Delta-1)^{-1}\partial_i
=
\partial_i(\Delta-1)^{-1}.\]
These two properties will play a key role in the proof that follows.

In addition, we will need the following claim that follows from the product rule for derivatives.

\begin{thm}{Claim}\label{ex:hoelder-norm}
Given a function $v\in C^\infty(\TT^n)$, there is a sequence of constants $a_0,a_1,\dots$ such that
\[\|u\cdot v\|_{C^{k,\alpha}}\le a_0\cdot\|u\|_{C^{k,\alpha}}+a_k\cdot\|u\|_{C^{k-1,\alpha}}\]
for any integer $k\ge 1$ and $u\in C^\infty(\TT^n)$.
\end{thm}

\paragraph{Proof of \ref{lem:nash-guenter:guenter}.}
Recall that
$Q_{ij}(w,w)=\langle\partial_i w,\partial_jw\rangle$.
Let us show that 
\[(\Delta-1) Q_{ij}(w,w)=A_{ij}(w)
+\partial_iA_j(w)
+\partial_j A_i(w),\eqlbl{eq:Delta}\]
where each of $A_{ij}$, $A_{i}$, and $A_{j}$ is a linear combination (with constant coefficients) of the following quadratic forms
\[
\begin{aligned}
w&\mapsto \langle\partial_aw,\partial_bw\rangle,
&
w&\mapsto \langle\partial_{ab}w,\partial_cw\rangle,
&
w&\mapsto \langle\partial_{ab}w,\partial_{cd}w\rangle.
\end{aligned}
\eqlbl{eq:forms}
\]
Here $a,b,c,d$ range over all indices.

Indeed,
\begin{align*}
(\Delta-1) Q_{ij}(w,w)&=
\langle\partial_{i}\Delta w,\partial_jw\rangle
+\langle\partial_{i} w,\partial_{j}\Delta w\rangle+
\\
&\qquad\qquad+2\cdot \sum_s\langle\partial_{is} w,\partial_{js}w\rangle
-\langle\partial_{i} w,\partial_{j}w\rangle
=
\\
&=
2\cdot \sum_s
\langle\partial_{is} w,\partial_{js}w\rangle
-\langle\partial_{i} w,\partial_{j}w\rangle
-2\cdot \langle\Delta w,\partial_{ij}w\rangle
+
\\
&\qquad\qquad+\partial_i\langle\Delta w,\partial_{j}w\rangle
+\partial_j\langle\Delta w,\partial_{i} w\rangle.
\end{align*}
This proves \ref{eq:Delta}.

Now choose $A=A_i$ or $A=A_{ij}$ for some $i$ and $j$.
There is a sequence of constants $c_2,c_3,\dots$ such that
\[
\begin{aligned}
\|A(w)\|_{C^{0,\alpha}}&\le c_2\cdot \|w\|^2_{C^{2,\alpha}},
\\
\|A(w)\|_{C^{k-2,\alpha}}&\le c_2\cdot \|w\|_{C^{k,\alpha}}\cdot \|w\|_{C^{2,\alpha}}+c_{k}\cdot\|w\|^2_{C^{k-1,\alpha}}.
\end{aligned}
\eqlbl{eq:A=<C}\]
These inequalities follow by checking each form in \ref{eq:forms}.
Let us do the case $w\mapsto \langle\partial_{ab}w,\partial_{cd}w\rangle$.
The first inequality in \ref{eq:A=<C} is trivial.
The second follows since for any partial derivative $\partial^{\mathcal{I}}$ of order $k-2$ we have
\[\partial^{\mathcal{I}}\langle\partial_{ab}w,\partial_{cd}w\rangle=\langle\partial_{ab}\partial^{\mathcal{I}}w,\partial_{cd}w\rangle+\langle\partial_{ab}w,\partial_{cd}\partial^{\mathcal{I}}w\rangle + \text{extra terms},\]
where each extra term has the form $\langle\partial^{\mathcal{K}}w,\partial^{\mathcal{L}}w\rangle$ and partial derivatives
$\partial^{\mathcal{K}}$ and $\partial^{\mathcal{L}}$ have order at most $k-1$.

By \ref{prop:Delta-1}, we have another sequence of constants $c_2, c_3, \dots$ such that
\[
\begin{aligned}
\|(\Delta-1)^{-1}A(w)\|_{C^{2,\alpha}}&\le c_2\cdot \|w\|_{C^{2,\alpha}}^2,
\\
\|(\Delta-1)^{-1}A(w)\|_{C^{k,\alpha}}&\le c_2\cdot \|w\|_{C^{k,\alpha}}\cdot \|w\|_{C^{2,\alpha}} +c_{k}\cdot\|w\|_{C^{k-1,\alpha}}^2.
\end{aligned}
\eqlbl{eq:AB=<C}\]

Recall also that
\[
L_{ij}(v)=-2\cdot\langle\partial_{ij} f,v\rangle
+\partial_{i}\langle\partial_{j} f,v\rangle
+\partial_{j}\langle\partial_{i} f,v\rangle;
\eqlbl{eq:L-def2}
\]
see the proof of Nash's lemma.
Let us define
$\tilde Q(w,w)$ as a linear combination of $\partial_if$ and $\partial_{ij}f$ for all $i\le j$
such that
\[
\begin{aligned}
-2\cdot\langle \partial_{ij}f,\tilde Q(w,w)\rangle
&=(\Delta-1)^{-1} A_{ij}(w),
\\
\langle \partial_if,\tilde Q(w,w)\rangle
&=(\Delta-1)^{-1} A_i(w)
\end{aligned}
\eqlbl{eq:tildeQ}
\]
for all $i\le j$.
Since $f$ is free, $\tilde Q(w,w)$ is well-defined, and it extends to a symmetric bilinear form.
Combining \ref{eq:L-def2}, \ref{eq:tildeQ}, and \ref{eq:Delta}, we get
\[
L\tilde Q=Q.
\] 

Let $\{e_i,e_{ij}\}_{i\le j}$ be the dual frame to $\{\partial_if,\partial_{ij}f\}_{i\le j}$ in the osculating bundle $\T^2$.
Notice that
\[\tilde Q(w,w)=\sum_i(\Delta-1)^{-1} A_i(w)\cdot e_i-\tfrac12\cdot\sum_{i\le j}(\Delta-1)^{-1}A_{ij}(w)\cdot e_{ij}.\eqlbl{eq:tildeQend}\]
Since the maps $e_i, e_{ij}\colon \TT^n \to \RR^d$ are smooth, \ref{eq:AB=<C} together with \ref{ex:hoelder-norm} imply the inequalities in \ref{lem:nash-guenter:guenter};
in particular, the coefficient in front of $\|w\|_{C^{2,\alpha}}\cdot\|w\|_{C^{k,\alpha}}$ remains unchanged for all $k$.
\qeds

The promised tangential term $\vec x=\vec x(w)$ can be defined as a linear combination of $\partial_if$ such that
\[\langle \partial_if,\vec x\rangle=-(\Delta-1)^{-1} A_i(w)\]
for each $i$; compare with the second line in \ref{eq:tildeQ}.

\paragraph{Acknowledgments.}
This work was partially supported by the National Science Foundation, grant DMS-2005279.

I would like to thank
the anonymous referee,
Mikhael Gromov,
Fedor Nazarov,
and Deane Yang for their help.

{\small\sloppy
\def\emph{\textit}

\printbibliography[heading=bibintoc]
\fussy
}

\end{document}